\documentclass{agtart_a}
\pdfoutput=1

\usepackage{graphicx}
\usepackage{mathdots}
\usepackage{color}

\title{One-sided Heegaard splittings of $\mathbb{R}\mathrm{P}^3$}

\author{Loretta Bartolini}
\givenname{Loretta}
\surname{Bartolini}
\address{Department of Mathematics and Statistics\\
University of Melbourne\\\newline
Parkville VIC 3010\\Australia}
\email{L.Bartolini@ms.unimelb.edu.au} 
\urladdr{}

\author{J\,Hyam Rubinstein}
\givenname{J\,Hyam}
\surname{Rubinstein}
\email{H.Rubinstein@ms.unimelb.edu.au}
\urladdr{}

\volumenumber{6}
\issuenumber{}
\publicationyear{2006}
\papernumber{46}
\startpage{1319}
\endpage{1330}

\doi{}
\MR{}
\Zbl{}

\keyword{one-sided Heegaard splitting}
\keyword{geometrically compressible}
\subject{primary}{msc2000}{57M27}
\subject{secondary}{msc2000}{57N10}

\received{31 August 2005}
\revised{6 July 2006}
\accepted{20 July 2006}
\published{20 September 2006}
\publishedonline{20 September 2006}
\proposed{}
\seconded{}
\corresponding{}
\editor{}
\version{}

\arxivreference{math.GT/0509007}

\AtBeginDocument{\let\bar\wbar\let\tilde\wtilde\let\hat\what}
\makeop{genus}

\makeatletter
\def\cnewtheorem#1[#2]#3{\newtheorem{#1}{#3}[section]
\expandafter\let\csname c@#1\endcsname\c@thm}
\makeatother

\newcounter{step}
\newcounter{claim}
\newcounter{subclaim}

\newtheorem{thm}{Theorem}[section]
\cnewtheorem{corollary}[thm]{Corollary}
\cnewtheorem{proposition}[thm]{Proposition}
\cnewtheorem{lem}[thm]{Lemma}
\theoremstyle{definition}
\cnewtheorem{defn}[thm]{Definition}

\newenvironment{step}{\medskip \textbf{Step \arabic{step}}\qua \stepcounter{step}}{\medskip}
\newenvironment{claim}{\medskip \textbf{Claim}\qua}{\medskip}
\newenvironment{subclaim}{\medskip \textbf{Subclaim \arabic{subclaim}}\qua \stepcounter{subclaim}}{\medskip}

\newcommand{\hs}{Heegaard splitting}
\newcommand{\oshs}{one-sided Heegaard splitting}

\newcommand{\oss}{one-sided splitting}
\newcommand{\tshs}{two-sided Heegaard splitting}
\newcommand{\tss}{two-sided splitting}
\newcommand{\gi}{geometrically incompressible}
\newcommand{\gc}{geometrically compressible}

\newcommand{\rpa}{$\mathbb{R}\mathrm{P}^2$}
\newcommand{\rpb}{$\mathbb{R}\mathrm{P}^3$}

\begin{document}

\begin{asciiabstract}
Using basic properties of one-sided Heegaard splittings, a direct
proof that geometrically compressible one-sided splittings of RP^3 are
stabilised is given. The argument is modelled on that used by
Waldhausen to show that two-sided splittings of S^3 are standard.
\end{asciiabstract}

\begin{htmlabstract}
Using basic properties of one-sided Heegaard splittings, a direct
proof that geometrically compressible one-sided splittings of
<b>R</b>P<sup>3</sup> are stabilised is given. The argument is modelled on
that used by Waldhausen to show that two-sided splittings of S<sup>3</sup> are
standard.
\end{htmlabstract}

\begin{abstract} 
Using basic properties of one-sided Heegaard splittings, a direct
proof that geometrically compressible one-sided splittings of
$\mathbb{R}P^3$ are stabilised is given. The argument is modelled on
that used by Waldhausen to show that two-sided splittings of $S^3$ are
standard.
\end{abstract}

\maketitle

\setcounter{step}{1}
\setcounter{subclaim}{1}

\section{Introduction}

Since their formal introduction in 1978~\cite{rubin78}, \oshs{s} of $3$--mani\-folds have been the subject of little study. This paucity of literature can largely be attributed to the lack of generality of such splittings, as compared with classical \hs{s}, and the invalidity of an analogue to Dehn's lemma and the loop theorem~\cite{stallings60}. Various works, both prior and subsequent to Rubinstein \cite{rubin78}, have addressed non-orientable surfaces in $3$--manifolds such as Bredon and Wood \cite{bredon-wood}, Hempel \cite{hempel74}, Frohman \cite{frohman86} and Rannard \cite{rannard96}, and classifications are made in the latter works when restricted to \gi{} surfaces. However, in order to study \oss{s} effectively, the existence and behaviour of \gc{} splittings must be considered. 

Well known in \tshs{} theory, the stabilisation problem is also present for \oss{s}. By its very nature, this issue demands an understanding of \gc{} splitting surfaces. To date, no connection has been drawn between geometric compressibility and stabilisation. Here, a direct correspondence is drawn for the simplest case: \rpb{}. 

The result is analogous to that of Waldhausen's for \tss{s} of $S^3$~\cite{wald68} and it is upon these original arguments that the proof is based. While there have been many subsequent proofs of the $S^3$ case using simpler arguments, in the absence of an analogue to Casson and Gordon's result on weak reducibility, such approaches are not currently viable for \oss{s}.

We would like to thank Marc Lackenby for helpful discussions and feedback that assisted in the preparation of this paper.

\section{One-sided Heegaard splittings}
Throughout, let $M$ be a closed, orientable $3$--manifold and consider all manifolds and maps as PL.

\begin{defn}
A pair $(M, K)$ is called a \oshs{} if $K$ is a closed non-orientable surface embedded in $M$ such that $H = M \setminus K$ is an open handlebody.
\end{defn}

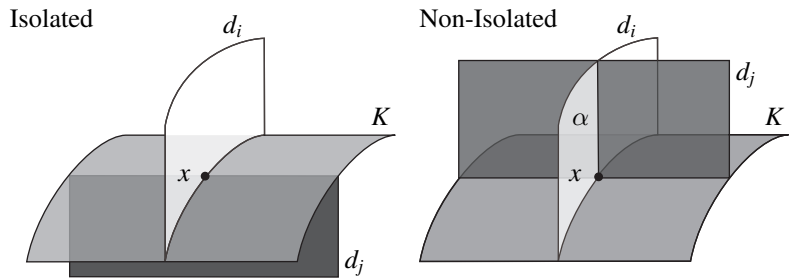
\begin{figure}[ht!]\small
   \centering
   \includegraphics[width=4in]{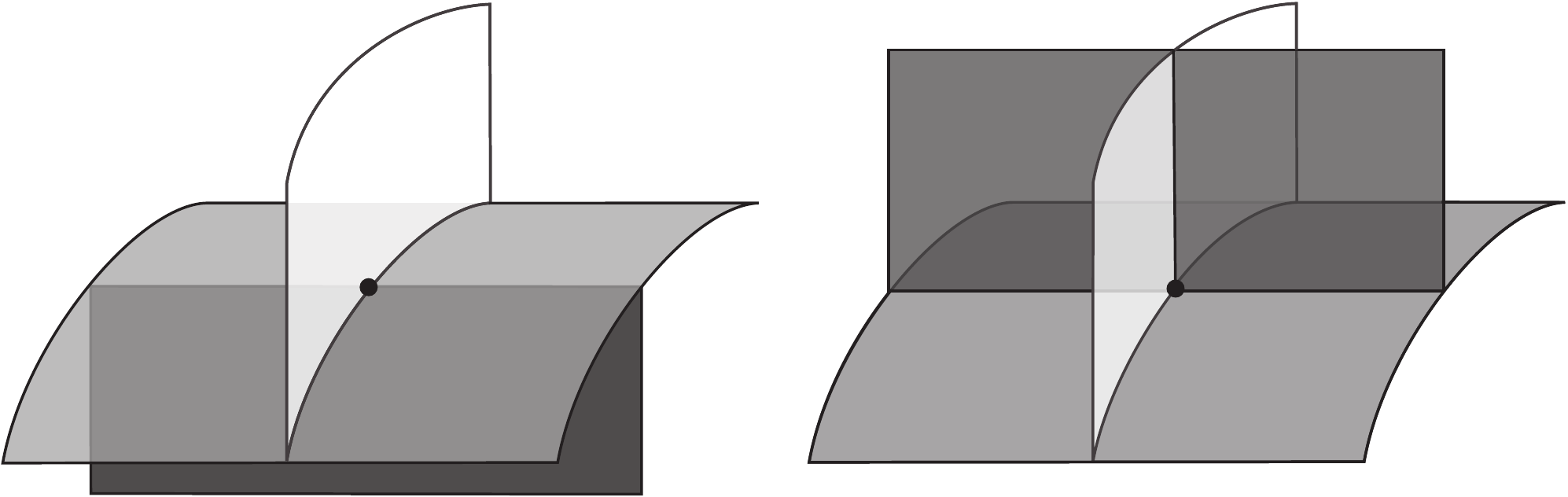} 
   \put(-295, 95){Isolated}
   \put(-140, 95){Non-Isolated}
   \put(-215, 93){$d_i$}
   \put(-66, 93){$d_i$}
   \put(-21, 75){$d_j$}
   \put(-159, 58){$K$}
   \put(-10, 58){$K$}
   \put(-83, 58){$\alpha$}
   \put(-232, 37){$x$}
   \put(-84, 37){$x$}
   \put(-169, 4){$d_j$}
   \caption{Different intersection types for meridian discs of $(M, K)$}
   \label{fig:intersect}
\end{figure}

As with \tss{s}, it is useful to consider {\it meridian discs\/} for $(M, K)$, which are taken to be the closure of meridian discs for the handlebody complement $H$ in the usual sense. Due to the non-orientability of $K$, the boundaries of such discs can intersect themselves or one another in two distinct ways (see \fullref{fig:intersect}).

\begin{defn} If $x = \partial d_i \cap \partial d_j$, where $d_i, d_j$ are meridian discs for a \oss{}, and $B_{\varepsilon}(x)$ is a small ball centred at $x$, call $x$ isolated if $d_i \cap d_j \cap B_{\varepsilon}(x) = x$. Call $x$ non-isolated if $d_i \cap d_j \cap B_{\varepsilon}(x) = \alpha$, where $\alpha$ is an arc containing $x$.
\end{defn}

\subsection{Existence}

\begin{thm}[Rubinstein \cite{rubin78}]\label{thm:existence}
For any element $\alpha \not= 0$ in $H_2(M, \mathbb{Z}_2)$, there is a \oshs{} $(M, K)$ with $[K] = \alpha$.
\end{thm}

The \oss{} technique is hence applicable to a large class of $3$--manifolds, which can be easily identified using algebraic methods. Associated with any \oss{} is a double cover $p\co \tilde{M} \rightarrow M$, where $\tilde{K} = p^{-1}(K)$ is the orientable double cover of $K$. The surface $\tilde{K}$ gives a natural \tss{} of $\tilde{M} = p^{-1}(M)$, with handlebody components interchanged by the covering translation $g\co \tilde{M} \rightarrow \tilde{M}$. 

In order to consider the simplest surface representing a $\mathbb{Z}_2$--homology class, a notion of incompressibility for non-orientable surfaces is required.

\begin{defn}
A surface $K \not= S^2$ embedded in $M$ is \textit{\gi{}\/} if any simple, closed, noncontractible loop on $K$ does not bound an embedded disc in $M$. Call $K$ \gc{} if it is not \gi{}.
\end{defn}

The existence of such a \oss{} surface is not implied by existence of \oss{s} in general.  However, by restricting to the class of irreducible, non-Haken $3$--manifolds, such a connection can be drawn.

\begin{thm}[Rubinstein \cite{rubin78}]\label{thm:incomp}
If $M$ is irreducible and non-Haken, then there is a \gi{} \oss{} associated with any nonzero class in $H_2(M, \mathbb{Z}_2)$.
\end{thm}

While little is known about general \gi{} one-sided surfaces in $3$--manifolds, a classification is available for Seifert fibered spaces. The Lens space case is discussed by the second author \cite{rubin78} and general Seifert fibered spaces in Frohman \cite{frohman86} and Rannard \cite{rannard96}. Considering \rpb{} as $L(2, 1)$, the former result is sufficient here.

Combining \fullref{thm:existence} and \fullref{thm:incomp}, any Lens space of the form $L(2k, q)$, where $(2k, q) = 1$, has \gi{} \oshs{s}. In~\cite{rubin78}, it is shown that any such space has a unique, \gi{} splitting that realises the minimal genus of all \oss{s} of the manifold. An algorithm is given by Bredon and Wood \cite{bredon-wood} for calculating this genus. Since $H_2(L(2k, q); \mathbb{Z}) = 0$ and all \oss{} surfaces of a Lens space are represented by the same $\mathbb{Z}_2$--homology class, any splitting surface that is \gc{} must geometrically compress to the minimal genus surface.

\subsection{Stabilisation}

\begin{defn}
A \oss{} $(M, K)$ is stabilised if and only if there exists a pair of embedded meridian discs $d, d^{\prime}$ for $H$ such that $d \cap d^{\prime}$ is a single isolated point.
\end{defn}

\begin{defn}
A \oss{} of an irreducible manifold is called irreducible if it is not stabilised.
\end{defn}

As stabilised \oss{} surfaces are inherently \gc{}, irreducibility is implied by geometric incompressibility. In future work, we hope to give evidence that geometric incompressibility of \oss{} surfaces is actually analogous to strong irreducibility in the two-sided case.

\subsection{Stable equivalence}

\begin{defn}
One-sided Heegaard splittings $(M_1, K_1)$ and $(M_2, K_2)$ are equivalent if there exists a homeomorphism from $M_1$ to $M_2$ that maps $K_1$ to $K_2$.
\end{defn}

As for \tss{s}, there is a notion of stabilising distinct \oss{s} until they are equivalent. Let $(S^3, L)$ denote the standard genus 1 \tss{} of the $3$--sphere and $(M, K)\#n(S^3, L)$ be the connected sum of $(M, K)$ with $n$ copies of $(S^3, L)$.

\begin{defn}
One-sided splittings $(M_1, K_1)$ and $(M_2, K_2)$ are stably equivalent if $(M_1, K_1)\#n(S^3, L)$ is equivalent to $(M_2, K_2)\#m(S^3, L)$ for some $m, n$.
\end{defn}

Unlike \tss{s}, stable equivalence does not hold for \oshs{s} in general.  However, a version applies to splitting surfaces represented by the same $\mathbb{Z}_2$--homology class:

\begin{thm}[Rubinstein \cite{rubin78}]
If $(M, K_1)$ and $(M, K_2)$ are \oshs{s} with $[K_1] = [K_2]$, then they are stably equivalent.
\end{thm}

Motivated by the fact that the little that is known about \oshs{s} is largely restricted to \gi{} splitting surfaces, we use these basic properties of \oss{s} to broach geometric compressibility. Given any stabilised \oss{} is inherently \gc{}, it is natural to ask when geometric compressibility corresponds to stabilisation.

\section[One-sided Heegaard splittings of RP3]{One-sided Heegaard splittings of \rpb{}}

Investigating any existence of a correlation between geometric compressibility and stabilisation, the simplest case to consider is \rpb{}, which corresponds to $S^3$ in the two-sided case. Here, the original arguments given by Waldhausen are adapted to show that all \gc{} splittings of \rpb{} are stabilised. 

In brief, the approach is to take an unknown splitting and the known minimal genus splitting by \rpa{} and stabilise the two until they are equivalent. Keeping track of the disc systems introduced by this process, it is possible to arrange them such that the reverse process of destabilising to get the unknown splitting preserves dual pairs from the minimal genus splitting. Thus, dual discs exist for the original unknown splitting, hence it is stabilised.

\begin{thm}
Every \gc{} \oshs{} of \rpb{} is stabilised.
\end{thm}

\begin{proof}
Take a \gc{} \oshs{} $(M, K)$ of $M \cong$ \rpb{} and let $(M, P)$ be the splitting along $P \cong$ \rpa{}. Since $H_2(M; \mathbb{Z}_2) \cong \mathbb{Z}_2$, there is only one nontrivial $\mathbb{Z}_2$--homology class so $[K] = [P]$. As $P$ is the unique \gi{} splitting surface of $M$, the unknown splitting surface $K$ geometrically compresses to $P$.

By stable equivalence, each splitting surface can be stabilised a finite number of times until the two are equivalent. Represent this splitting by $(M, K^{\prime})$ and let $H = M \setminus K^{\prime}$ be the handlebody complement. Let $\Delta^{K}$ be a set of meridian discs introduced by stabilisations of $(M, K)$, chosen such that $\Delta^{K} = \Delta_K \cup \Delta_{\smash{K}}^{\prime}$, where $\Delta_K = d_1, d_2, \ldots, d_k$ and $\Delta_{\smash{K}}^{\prime} = d_1^{\prime}, d_2^{\prime}, \ldots, d_{\smash{k}}^{\prime}$ are sets of disjoint discs with $|d_i \cap d_{\smash{i}}^{\prime}| =   1$ and $d_i \cap d_{\smash{j}}^{\prime} = \emptyset$ for $i \not= j$. Then $|\Delta^K| = 2k = (\genus(K^{\prime}) - \genus(K))$. Note that this number is always even, as each stabilisation increases the genus of the handlebody by $2$.

Similarly, let $\Delta^P = \Delta_P \cup \Delta_{\smash{P}}^{\prime}$ be the set of discs introduced by stabilising $(M, P)$. Notice that since $M \setminus P$ is an open $3$--cell, $\Delta^P$ is a complete disc system for $H$.

Consider the non-isolated intersections between discs in $\Delta_K$, $\Delta_K^{\prime}$ and $\Delta_P$, $\Delta_P^{\prime}$. Let
$$\Lambda_0 = \{d \cap D\},\quad\Lambda_0^{\prime} = \{d^{\prime} \cap D^{\prime}\},
\quad\Lambda_1 = \{d \cap D^{\prime}\}\quad\text{and}\quad\Lambda_1^{\prime} = \{d^{\prime} \cap D\}$$
be the collections of arcs of intersection between the given pairs for all $d \in \Delta_K$, $d^{\prime} \in \Delta_K^{\prime}$, $D \in \Delta_P$, $D^{\prime} \in \Delta_P^{\prime}$.

Stabilise $(M, K^{\prime})$ along $\Lambda_0$, $\Lambda_0^{\prime}$, $\Lambda_1$, $\Lambda_1^{\prime}$. Call the resulting splitting $(M, K^{\prime\prime})$, with handlebody complement $H^{\prime} = M \setminus K^{\prime\prime}$. Let
\begin{center}
\begin{tabular}{ccccccccc}

$\bar{\Delta}_K$ && $\Delta_K$&&$\Lambda_0$, $\Lambda_1$ & & $\Delta_K^{\prime}$ & & $\Lambda_0^{\prime}$, $\Lambda_1^{\prime}$\\
$\bar{\Delta}_P$ &be& $\Delta_P$ & cut & $\Lambda_0$, $\Lambda_1^{\prime}$ & \textit{plus\/} the discs & $\Delta_P^{\prime}$&  along & $\Lambda_0^{\prime}$, $\Lambda_1$\\
$\bar{\Delta}_K^{\prime}$ && $\Delta_K^{\prime}$ &along& $\Lambda_0^{\prime}$, $\Lambda_1^{\prime}$&dual to cuts of& $\Delta_K$ && $\Lambda_0$, $\Lambda_1$\\
$\bar{\Delta}_P^{\prime}$ && $ \Delta_P^{\prime}$ && $\Lambda_0^{\prime}$, $\Lambda_1$&& $\Delta_P$ && $\Lambda_0$, $\Lambda_1^{\prime}$\\

\end{tabular}
\end{center}
where a disc dual to a cut along an arc $\lambda$ is a transverse cross-section of a closed regular neighbourhood of $\lambda$ (see \fullref{fig:stabilisation}). For such discs, use parallel copies for the $K$ and $P$ systems in order to retain dual pairs in each. Let $\bar{\Delta}^K = \bar{\Delta}_K \cup \bar{\Delta}_{\smash{K}}^{\prime}$ and $\bar{\Delta}^P = \bar{\Delta}_P \cup \bar{\Delta}_P^{\prime}$. Notice that $\bar{\Delta}^P$ is again a complete disc system for $H^{\prime}$.

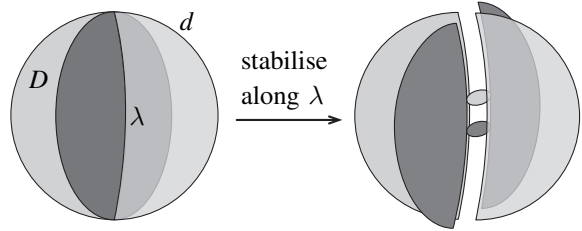
\begin{figure}[ht!]\small
   \centering
   \includegraphics[width=3in]{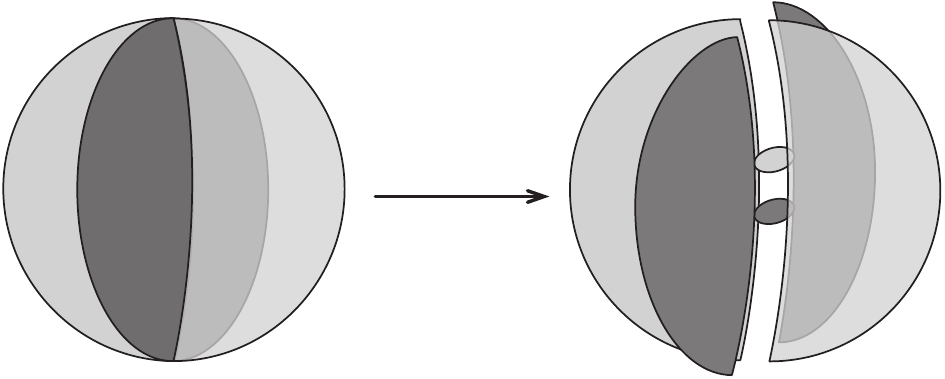}
   \put(-153, 76){$d$}
   \put(-210, 54){$D$}
   \put(-129, 61){stabilise}
   \put(-129, 48){along $\lambda$}
   \put(-172, 40){$\lambda$}
   \caption{Stabilising along an arc $\lambda$, where $d \in \Delta_K$ and $D \in \Delta_P$ or $\Delta_P^{\prime}$}
   \label{fig:stabilisation}
\end{figure}

The aim of this second stabilisation process is to remove all existing non-isolated intersections between $\Delta^P$ and $\Delta^K$. Therefore, it is imperative that the disc systems are not moved once this second set of stabilisations is complete, as any moves may introduce new intersections. Hence, the standard procedure of manipulating stabilising discs to get sets of disjoint dual pairs is not performed.

Order the $\bar{\Delta}_K, \bar{\Delta}_{\smash{K}}^{\prime}$ and $\bar{\Delta}_P, \bar{\Delta}_{\smash{P}}^{\prime}$ disc systems with respect to the nesting of arcs of stabilisation. For example, consider $\bar{d}_i, \bar{d}_j \in \bar{\Delta}_K$ that were split off $d \in \Delta_K$ by arcs $\lambda_i, \lambda_j$ respectively. If $\lambda_i$ is outermost with respect to the point $d \cap d^{\prime}$, then $j < i$ (see \fullref{fig:labelling}). Note that there is a rooted tree dual to the subdisc system for $d$, where the point of $d \cap d^{\prime}$ is the root, which induces the ordering. Label the dual discs such that $\bar{d}_k^{\prime}$ is a transverse cross-section of $\lambda_k$, hence $\bar{d}_k^{\prime} \in \bar{\Delta}_{\smash{K}}^{\prime}$ is dual to $\bar{d}_k$. Apply similar labelling to the $\bar{\Delta}_P, \bar{\Delta}_{\smash{P}}^{\prime}$ systems.

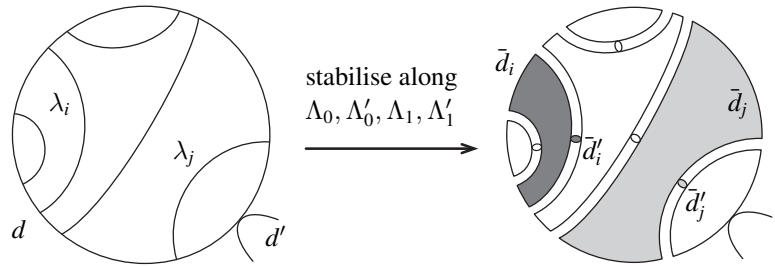
\begin{figure}[ht!]\small
   \centering
   \includegraphics[width=4in]{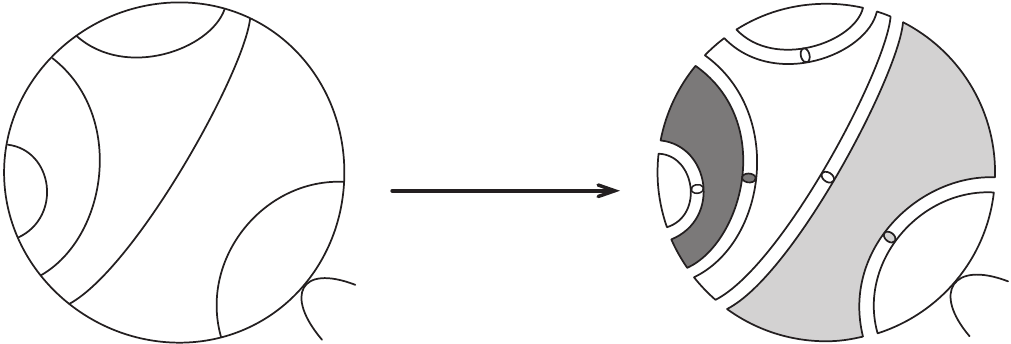} 
   \put(-106, 73){$\bar{d}_i$}
   \put(-177, 68){stabilise along}
   \put(-178, 55){$\Lambda_0, \Lambda_0^{\prime}, \Lambda_1, \Lambda_1^{\prime}$}
   \put(-275, 58){$\lambda_i$}
   \put(-18, 58){$\bar{d}_j$}
   \put(-228, 40){$\lambda_j$}
   \put(-73, 40){$\bar{d}_i^{\prime}$}
   \put(-34, 20){$\bar{d}_j^{\prime}$}
   \put(-289, 10){$d$}
   \put(-193, 7){$d^{\prime}$}
   \caption{Discs $\bar{d}_i, \bar{d}_j$ obtained by splitting $d$ along $\lambda_i, \lambda_j$, where $j < i$}
   \label{fig:labelling}
\end{figure}

Consider the intersections between discs $\bar{d}_i \in \bar{\Delta}_K$ and $\bar{d}_{\smash{j}}^{\prime} \in \bar{\Delta}_{\smash{K}}^{\prime}$. By construction, $\partial \bar{d}_i \cap \partial \bar{d}_{\smash{i}}^{\prime}$ is a single isolated point and $\partial \bar{d}_i \cap \{\partial \bar{d}_{\smash{j}}^{\prime}\ |\ j = 1, 2, ... , (i-1)\} = \emptyset$. For $i \leq j$, points of $\partial \bar{d}_i \cap \partial \bar{d}_j^{\prime}$ are isolated.

If $m = |\bar{\Delta}_K| = |\bar{\Delta}_{\smash{K}}^{\prime}|$, then $2m$ is the total change in genus from $K$ to $K^{\prime\prime}$. Construct the $2m \times 2m$ intersection matrix $\mathbf{M} = [m_{ij}]$ for discs in $\bar{\Delta}^K$. Define $m_{ij}$ as follows, where $|\partial \bar{d}_i \cap \partial \bar{d}_i|$ is given to be the number of isolated singularities of $\bar{d}_i$:
\begin{displaymath}
m_{ij} =
\left\{ \begin{array}{ll}
|\partial \bar{d}_i \cap \partial \bar{d}_j^{\prime}|, &1 \leq i, j \leq m \\
|\partial \bar{d}_{i} \cap \partial \bar{d}_{j-m}|,& 1 \leq i \leq m, (m+1) \leq j \leq 2m \\
|\partial \bar{d}_{i-m}^{\prime} \cap \partial \bar{d}_j^{\prime}|,& (m+1) \leq i \leq 2m, 1 \leq j \leq m \\
|\partial \bar{d}_{i-m}^{{\prime}} \cap \partial \bar{d}_{j-m}|,& (m+1) \leq i, j \leq 2m \\
\end{array} \right.
\end{displaymath}
Since $\bar{\Delta}_K$, $\bar{\Delta}_{\smash{K}}^{\prime}$ are systems of embedded, disjoint discs, the off-diagonal blocks are zero. By symmetry, the diagonal blocks are mutually transpose. While initially this symmetry makes the full matrix unnecessary, the asymmetry of later moves requires the consideration of all entries as described. Given the discs are not to be manipulated after the second set of stabilisations, the matrix is the identity if and only if $\{\Lambda_0, \Lambda_0^{\prime}, \Lambda_1, \Lambda_1^{\prime}\} = \emptyset$.
Thus:
\begin{displaymath}
\mathbf{M} = 
\left( \begin{array}{ccccc|ccccc}
1      & \star  & \star  & \ldots & \star  & 0 &\ldots&\ldots&\ldots&0 \\
0      & 1      & \star  & \ldots & \star  &&  &&& \\
\vdots & \ddots & \ddots & \ddots & \vdots &\vdots&& \iddots &&\vdots \\
0      & \ldots & 0      & 1      & \star  &&&&  & \\
0      & 0      & \ldots & 0      & 1      &0&\ldots&\ldots&\ldots& 0 \\
\hline
0 &\ldots&\ldots&\ldots&0&      1      &  0     & \ldots & 0      & 0      \\
&&&&&     \star  & 1      & 0      & \ldots & 0      \\
\vdots&& \iddots&&\vdots& \vdots & \ddots & \ddots & \ddots & \vdots \\
&&&&&     \star  & \ldots & \star  & 1      & 0      \\
0&\ldots&\ldots&\ldots& 0 &     \star  & \ldots & \star  & \star  & 1      \\
\end{array} \right)
\end{displaymath}
If $n = |\bar{\Delta}_P| = |\bar{\Delta}_{\smash{P}}^{\prime}|$, the $2n \times 2n$ intersection matrix $\mathbf{N}$ for the discs in $\bar{\Delta}^P$ can be constructed similarly. This $\mathbf{N}$ has a similar block structure to $\mathbf{M}$. 

Let $D = \bar{D}_n \in \bar{\Delta}_P$, the disc corresponding to the last row of the upper half of $\mathbf{N}$, and let $D^{\prime} \in \bar{\Delta}_{\smash{P}}^{\prime}$ be its dual. Thus $D, D^{\prime}$ are a dual pair disjoint from all other discs in $\bar{\Delta}^P$. However, several possibilities exist for how $D, D^{\prime}$ may intersect $\bar{\Delta}^K$:

\begin{enumerate}
\renewcommand{\labelenumi}{(\alph{enumi})}
\item Both $D, D^{\prime}$ are disjoint from $\bar{\Delta}^K$ or the pair intersect only one of $\bar{\Delta}_K, \bar{\Delta}_K^{\prime}$;
\item One of $D, D^{\prime}$ is disjoint from $\bar{\Delta}^K$, while the other intersects both $\bar{\Delta}_K, \bar{\Delta}_K^{\prime}$;
\item Both $D$ and $D^{\prime}$ intersect $\bar{\Delta}^K$ and Case (a) does not apply.
\end{enumerate}

In Case (a), compress along whichever of $\bar{\Delta}_K, \bar{\Delta}_K^{\prime}$ is disjoint from $D$ and $D^{\prime}$. This results in $(M, K)$, without having affected $D, D^{\prime}$, which remain a dual pair of embedded discs. Therefore, $(M, K)$ is stabilised.

In Case (b), suppose $D$ intersects both $\bar{\Delta}_K$ and $\bar{\Delta}_{\smash{K}}^{\prime}$, while $D^{\prime} \cap \bar{\Delta}^K = \emptyset$. Since $D^{\prime}$ is disjoint from $\bar{\Delta}_K$, it can be used to remove intersections between $D$ and $\bar{\Delta}_K$ by a process of band-summing:

Take $d \in \bar{\Delta}_K$, with $d \cap D \not= \emptyset$, such that there exists an arc $\alpha \subset \partial D$ with one endpoint at $D \cap D^{\prime}$, the other at $d \cap D$ and $\alpha \cap \bar{\Delta}_K = \emptyset$. Join a parallel copy of $D^{\prime}$ to $d$ by the boundary of a closed half-neighbourhood of $\alpha$. This removes one point from $d \cap D$ and since $\alpha \cap \bar{\Delta}_K = \emptyset$, no additional intersections are created within $\bar{\Delta}_K$. Repeat this procedure for all discs in $\bar{\Delta}_K$ that intersect $D$, taking care to work in an order that does not introduce intersections. Thus, all intersections between $D$ and $\bar{\Delta}_K$ can be removed without changing the intersection properties of $\bar{\Delta}_K$, resulting in Case (a) above.

In Case (c), both $D$ and $D^{\prime}$ intersect $\bar{\Delta}^K$, with no immediate means by which to remove intersections. Any attempts at band-summing, as used for Case (b), would introduce intersections between $\bar{\Delta}_K$ and $\bar{\Delta}_K^{\prime}$. Therefore, it is this case that requires significant attention.

\begin{claim}
After modifying $\bar{\Delta}^K$, there exists a dual pair of discs $\bar{d}, \bar{d}^{\prime} \in \bar{\Delta}^K$ such that $|\partial \bar{d} \cap \partial D| = 1$ and $|\partial \bar{d}^{\prime} \cap \partial D| \leq 1$ (or vice versa), and $D \cap (\bar{\Delta}^K \setminus \{\bar{d}, \bar{d}^{\prime}\}) = \emptyset$.
\end{claim}

The proof of this claim requires two steps, each of which is technical in nature. In particular, in the first step the most vital, yet most subtle, part of the argument appears.

\begin{step}
Describe surgery on $\bar{\Delta}^K$ in order to make $|\partial d \cap \partial D| \leq 1$ for all $d \in \bar{\Delta}^K$.
\end{step}

Consider arcs contained in $\partial D$ with endpoints on $\partial d$. Take a shortest arc $\alpha \subset \partial D$, with endpoints $\{a_0, a_1\}$ such that $a_i \in \partial d$ and $\alpha^{\circ} \cap d = \emptyset$. Such an arc can be chosen such that $\alpha \cap D^{\prime} = \emptyset$. If $\beta_1, \beta_2 \subset \partial d$ are the arcs with $\partial \beta_i = \{a_0, a_1\}$, let $\beta = \beta_i$ such that $\beta \cap d^{\prime}$ is a single point.

\begin{subclaim}
The loop $\gamma$ formed by the arcs $\alpha$ and $\beta$ bounds a disc in $H^{\prime}$ that is dual to $d^{\prime}$.
\end{subclaim}

In order to prove Subclaim 1, it is necessary to consider both isolated and non-isolated intersections between the $\bar{\Delta}^P$ and $\bar{\Delta}^K$ disc systems. For clarity, the subtleties are best captured by passing to the orientable double cover.

Take the orientable double cover $(\tilde{M}, \tilde{K}^{\prime\prime})$ corresponding to $(M, K^{\prime\prime})$, with covering projection $p\co \tilde{M} \rightarrow M$, covering translation $g\co \tilde{M} \rightarrow \tilde{M}$ and handlebody components $H_1, H_2$. Let $\tilde{d} = p^{-1}(d) \cap H_1$ and $\tilde{D} = p^{-1}(D) \cap H_2$, so an isolated intersection between $d$ and $D$ will correspond to discs in opposite handlebodies meeting in a point on the splitting surface. Let $\tilde{\beta} = p^{-1}(\beta) \cap \tilde{d}$ and $\tilde{\alpha} = p^{-1}(\alpha) \cap \tilde{D}$. Hence, the loop $\tilde{\gamma}$, bounded by $\tilde{\alpha}$ and $\tilde{\beta}$, is on $\tilde{K}^{\prime\prime}$ and constitutes part of the boundary of a disc in each handlebody.

Since all non-isolated intersections between $\bar{\Delta}^P$ and $\bar{\Delta}^K$ have been removed, the intersection $(p^{-1}(\bar{\Delta}^P) \cap H_i) \cap (p^{-1}(\bar{\Delta}^K) \cap H_i) = \emptyset$ for $i = 1$ or $2$. Specifically, $\tilde{d} \cap g(\tilde{D}) = g(\tilde{d}) \cap \tilde{D} = \emptyset$, so the loop $\gamma$ formed by $\alpha, \beta$ on $K^{\prime\prime}$ lifts to a pair of disjoint loops $\tilde{\gamma}, g(\tilde{\gamma})$ on $\tilde{K}^{\prime\prime}$ formed by $\tilde{\alpha}, \tilde{\beta}$ and $g(\tilde{\alpha}), g(\tilde{\beta})$ respectively.

As $\bar{\Delta}^P, \bar{\Delta}^K$ have no non-isolated intersections, $\tilde{d}$ is disjoint from $p^{-1}(\bar{\Delta}^P) \cap H_1$, which is a complete disc system for $H_1$. Thus the loop $\tilde{\gamma}$ bounds a disc $\tilde{d}^{\smash{1}}$ in $H_1$. Applying similar arguments to $g(\tilde{d})$ and $p^{-1}(\bar{\Delta}^P) \cap H_2$, the translated loop $g(\tilde{\gamma})$ bounds $g(\tilde{d}^1)$ in $H_2$. Since $\tilde{\gamma}, g(\tilde{\gamma})$ are disjoint, $\tilde{d}^1, g(\tilde{d}^1)$ are discs in opposite handlebodies with disjoint boundaries, hence $\tilde{d}^1 \cap g(\tilde{d}^1) = \emptyset$. Projecting to $(M, K^{\prime\prime})$, the disc $d^1 = p(\tilde{d}^1 \cup g(\tilde{d}^1))$ is embedded and dual to $d^{\prime}$, by choice of $\beta$.

If $\alpha$ is disjoint from $\bar{\Delta}^K \setminus d$, replace $d$ with $d^1$, which has two fewer points of intersection with $D$ than $d$. Repeat the process to remove all pairs of adjacent points in $d^1 \cap D$. Let $d_{\alpha}$ be the resulting disc and replace $d$ with $d_{\alpha}$ in $\bar{\Delta}_K$.

Any remaining arc $\alpha \subset \partial D$ between points of intersection with $d_{\alpha}$ is interrupted by intersections with $\bar{\Delta}^K \setminus d_{\alpha}$. These points of intersection are necessarily isolated.

\begin{subclaim}
There is a disc $d_0$ with $\partial \alpha_0 \subset \partial d_0$ for some $\alpha_0 \subset \alpha$ such that the intersection $\alpha_0 \cap (\bar{\Delta}^K \setminus d_0) = \emptyset$.
\end{subclaim}

Take $d_{\smash{K}} \in \bar{\Delta}^K \setminus d_{\alpha}$ with $x \in (d_K \cap \alpha)$ and again lift to the orientable double cover. Let $\tilde{d}_K = p^{-1}(d_K) \cap H_1$, so $p^{-1}(x) \in (\tilde{d}_K \cap \tilde{\alpha}) \cup (g(\tilde{d}) \cap g(\tilde{\alpha}))$ since $\tilde{d_K} \cap g(\tilde{D}) = \emptyset$. Now both $\tilde{d}_{\alpha}$ and $\tilde{d}_{\smash{K}}$ intersect $\tilde{D}$. By the previous argument, $\alpha$ and part of $\partial \tilde{d}_{\alpha}$ bound a disc $\tilde{d}_{\smash{\alpha}}^1$ in $H_1$. Since $\tilde{d}_{\smash{\alpha}}^1$ and $\tilde{d}_{\smash{K}}$ are discs in the same handlebody, their boundaries intersect in pairs of points. However, $\tilde{d}_K$ does not intersect $\tilde{d}_{\alpha}$, so both points of intersection lie on $\tilde{\alpha}$. Therefore, $\tilde{d}_{\smash{K}}$ intersects $\tilde{\alpha}$ in pairs of points. Similar arguments apply to $g(\tilde{d}_{\alpha}^1), g(\tilde{d}_K)$, thus $d_K$ intersects $\alpha$ in pairs of points.

Applying the above argument to any discs intersecting the subarc $\alpha_K \subset \alpha$, where $\partial \alpha_K \subset \partial d_K$, yields that any arcs of intersection between $\bar{\Delta}^K \setminus d_{\alpha}$ and $d_{\alpha}^1$ are nested. Therefore, there exists an innermost pair corresponding to intersections with the desired disc $d_0 \in \bar{\Delta}^K \setminus d_{\alpha}$ (see \fullref{fig:nesting}).

\begin{figure}[ht!]\small
   \centering
   \includegraphics[width=2in]{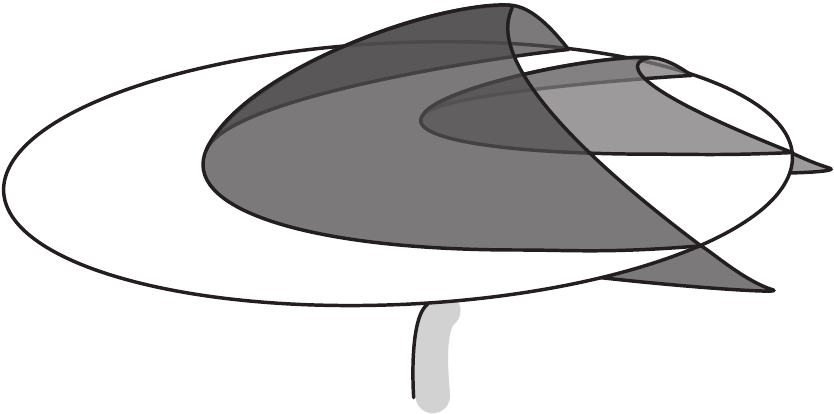}
   \put(-96, 69){$d_K$}
   \put(-27, 62){$d_0$}
   \put(-145, 57){$d_{\alpha}^1$}
   \put(-13, 55){$\alpha_0$}
   \put(-14, 30){$\alpha_K$}
   \put(-62, 13){$\alpha$}
   \put(-86, 3){$D$}
   \caption{Nested discs intersecting $\alpha$}
   \label{fig:nesting}
\end{figure}
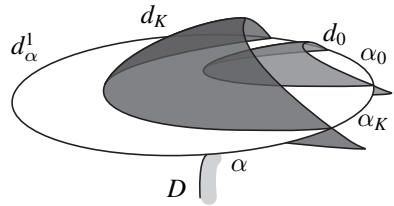

Apply the previous surgery to split $d_0$ along $\alpha_0$ and reduce the number of points of intersection with $D$. Continue this process, from edgemost arcs inwards, to remove all pairs of intersection points between $D$ and $\bar{\Delta}^K \setminus d_{\alpha}$. Applying the previous surgery to $d_{\alpha}$, the pair of intersection points that constitute the boundary of $\alpha$ can then be removed. Hence, the number of intersection points can be reduced to at most one.

Therefore, $D$ intersects any disc in $\bar{\Delta}^K$ in at most one point. If $D$ is disjoint from all such discs, then Case (a) above applies and the result holds.

\begin{step}
Reduce the number of discs in $\bar{\Delta}^K$ that have nonempty intersection with $D$ to at most a dual pair.
\end{step}

Consider two discs, each intersecting $D$ in a single point and let $\lambda \subset \partial D$ be the arc with an endpoint on each disc. Choose $d_a, d_b \in \bar{\Delta}^K$ to be such that $\lambda \cap (\bar{\Delta}^K \setminus \{d_a, d_b\}) = \emptyset$. Say that such discs are adjacent, since they are directly next to each other with respect to $\partial D$.

Do not perform surgery if the discs are a dual pair. Otherwise, take a parallel copy of whichever of $d_a, d_b$ corresponds to a later stabilisation---say $d_b$. Join the copy of $d_b$ to $d_{\smash{a}}$ by the boundary of a closed half-neighbourhood of $\lambda$. This forms a new disc $\bar{d}_a$ with $\bar{d}_a \cap D = \emptyset$. Replace $d_a$ in $\bar{\Delta}^K$ with $\bar{d}_a$. Note that any intersections of $d_b$ with $\bar{\Delta}^K$ will be present in $\bar{d}_a$.

The effect of the surgery on the intersection matrix is to add the row of $\mathbf{M}$ corresponding to $d_b$ to that corresponding to $d_a$. If both $d_a, d_b$ belong to one of $\bar{\Delta}_K$, $\bar{\Delta}_{\smash{K}}^{\prime}$, the surgery does not affect the off-diagonal blocks of $\mathbf{M}$. However, if $d_a = \bar{d}_k \in \bar{\Delta}_K$, $d_b = \bar{d}_l^{\prime} \in \bar{\Delta}_K^{\prime}$, where $k < l$, the $k$--th row of $\mathbf{M}$ becomes:
$$(\ \underbrace{0 \ \ \ldots \ \ 0}_{k} \ \ 1 \ \ \star \ \ \ldots\ldots \ \ \star \ | \ \star \ \ \ldots \ \ \star \ \ 1 \ \ \overbrace{0 \ \ \ldots \ \underbrace{0 \ \ldots \ \ 0}_{k}}^{l}\ )$$
Specifically, the $(m{+}k)$--th entry of the $k$--th row remains $0$. Therefore, throughout all surgery, discs in $\bar{\Delta}^K$ remain embedded. This allows the procedure to be iterated if necessary.

Perform surgery on all adjacent discs (except dual pairs) until there is, at most, a single pair of dual discs $\bar{d}, \bar{d}^{\prime}$, each intersecting $D$ in a single point. This pair corresponds to the latest stabilisation of any discs that had intersected $D$ after Step 1. Note that this may not be the pair corresponding to the centremost rows of $\mathbf{M}$, as these discs may not have initially intersected $D$. 

Having thus found discs $\bar{d}$ and $\bar{d}^{\prime}$ that prove the claim, it is now possible to destabilise $K^{\prime\prime}$ in a useful manner:

\begin{step} Replace $\bar{d}^{\prime}$ in $\bar{\Delta}_{\smash{K}}^{\prime}$ with $D$. Compress along $D$, thus destabilising $K^{\prime\prime}$. Discard $\bar{d}$.
\end{step}

Since $D$ is disjoint from $\bar{\Delta}^K \setminus \{\bar{d}, \bar{d}^{\prime}\}$, all other discs in this system remain intact after the compression. Therefore, the remaining discs again form systems of embedded dual pairs that correspond to stabilisations of $K$ and $P$, the latter of which is complete with respect to the newly destabilised splitting surface. As the original properties required for surgery on the discs systems are retained, Steps 1, 2 and 3 can be repeated for remaining discs in $\bar{\Delta}_P$. If the process is not terminated by the occurrence of Cases (a) or (b) as described previously, this process of destabilisation continues until it results in the original splitting $(M, K)$.

Since $(M, P)$ has minimal genus, $|\bar{\Delta}^K| < |\bar{\Delta}^P|$ as $K \not\cong$ \rpa{}. Therefore, after destabilising $(M, K^{\prime\prime})$ to get $(M, K)$ by the above process, there are dual pairs of discs remaining in $\bar{\Delta}^P$. Therefore, $(M, K)$ is stabilised. 
\end{proof}

\bibliographystyle{gtart}
\bibliography{link}

\end{document}